\documentclass[12pt,a4paper]{amsart}

\usepackage{anysize}
\usepackage[latin1]{inputenc}

\usepackage{hyperref}

\marginsize{3cm}{3cm}{3.5cm}{3.5cm}

\begin{document}
\title{A new Ramanujan-like series for $1/\pi^2$}

\author{Jesús Guillera}
\email{jguillera@gmail.com}
\address{Av. Cesáreo Alierta, 31 esc. izda {\rm $4^o$}--A, Zaragoza (Spain)}
\keywords{Hypergeometric series; WZ-method; Ramanujan-like series for $1/\pi^2$}
\date{}

\newtheorem{observacion}{Remark}

\begin{abstract}
Our main results are a WZ-proof of a new Ramanujan-like series for $1/\pi^2$ and a hypergeometric identity involving three series.
\end{abstract}

\maketitle

\section{The WZ-method}
We recall that a function $A(n,k)$ is hypergeometric in its two variables if the quotients
\[
\frac{A(n+1,k)}{A(n,k)} \quad {\rm and} \quad \frac{A(n,k+1)}{A(n,k)}
\]
are rational functions in $n$ and $k$, respectively. Also, a pair of hypergeometric functions in its two variables, $F(n,k)$ and $G(n,k)$ is said to be a Wilf and Zeilberger (WZ) pair \cite[Chapt. 7]{petkovsek} if
\begin{equation}\label{pro-WZ-pair}
F(n+1,k)-F(n,k)=G(n,k+1)-G(n,k).
\end{equation}
In this case, H. S. Wilf and D. Zeilberger \cite{wilf} have proved that there exists a rational function $C(n,k)$ such that
\begin{equation}\label{certificado}
G(n,k)=C(n,k)F(n,k).
\end{equation}
The rational function $C(n,k)$ is the so-called certificate of the pair $(F,G)$. To discover WZ-pairs, we use EKHAD \cite[Appendix A]{petkovsek}, a software written by D. Zeilberger. If EKHAD certifies a function, we have found a WZ-pair!. Then, if we sum (\ref{pro-WZ-pair}) over all $n \geq 0$, we get
\begin{equation}\label{sumas-wz}
\sum_{n=0}^{\infty} G(n,k) - \sum_{n=0}^{\infty} G(n,k+1) = -F(0,k) + \lim_{n \to \infty} F(n,k).
\end{equation}
We will write the functions $F(n,k)$ and $G(n,k)$ using rising factorials, also called Pochhammer symbols, rather than the ordinary factorials. The rising factorial is defined by
\begin{equation}\label{poch1}
(x)_n=\left\{
\begin{array}{ll}
x(x+1)\cdots(x+n-1), & \qquad n \in \mathbb{Z}^{+}, \\
1, & \qquad n=0,
\end{array} \right.
\end{equation}
or more generally by
\begin{equation}\label{poch2}
(x)_k=\frac{\Gamma(x+k)}{\Gamma(x)}.
\end{equation}
For $k \in \mathbb{Z}-\mathbb{Z}^{-}$, (\ref{poch2}) coincide with (\ref{poch1}). But (\ref{poch2}) is more general because it is also defined for all complex $x$ and $k$ such that
$x+k \in \mathbb{C}- ( \mathbb{Z}-\mathbb{Z}^{+} )$.
\par To use package EKHAD we will replace groups of rising factorials according to the following equivalences
\begin{align}
(1+k)_n &=\frac{(n+k)!}{k!}, \\
\left( \frac{1}{2}+k \right)_n &=\frac{1}{2^{2n}} \frac{(2n+2k)! k!}{(n+k)! (2k)!}, \\
\left( 1+\frac{k}{2} \right)_n \left( \frac{1}{2}+\frac{k}{2} \right)_n &=\frac{1}{2^{2n}} \frac{(2n+k)!}{k!}, \\
\left( 1+\frac{k}{3} \right)_n \left( \frac{1}{3}+\frac{k}{3} \right)_n \left( \frac{2}{3}+\frac{k}{3} \right)_n &=\frac{1}{3^{3n}} \frac{(3n+k)!}{k!},
\end{align}
which we can derive easily from the properties of the Gamma function.

\section{A new ramanujan-like series for $1/\pi^2$}
This paper is originated when we checked that EKHAD certifies the function
\begin{equation}
F(n,k)= \frac{\left( \frac{1}{2} \right)_n^3 \left( 1+\frac{k}{2} \right)_n \left( \frac{1}{2}+\frac{k}{2} \right)_n }{(1)_n^3 (1+k)_n^2} \frac{\left( \frac{1}{2} \right)_k^2}{(1)_k^2} \frac{96n^3}{2n+k},
\end{equation}
giving the companion
\begin{equation}
G(n,k)=\frac{\left( \frac{1}{2} \right)_n^3 \left(1+\frac{k}{2} \right)_n \left( \frac{1}{2}+\frac{k}{2} \right)_n }{(1)_n^3 (1+k)_n^2} \frac{\left( \frac{1}{2} \right)_k^2}{(1)_k^2} \, \frac{12k(8n^2+6kn+2n+k)}{2n+k}.
\end{equation}
As $F(0,k)=0$ and the last limit in (\ref{sumas-wz}) is also zero, we get
\begin{equation}\label{GigualG}
\sum_{n=0}^{\infty} G(n,k) = \sum_{n=0}^{\infty} G(n,k+1).
\end{equation}
As a consequence of Weierstrass M-test \cite[p. 49]{whittaker}, the convergence of this series is uniform. Therefore, the following steps hold
\begin{align}\label{uniform}
\lim_{k \to \infty} \sum_{n=0}^{\infty} G(n,k) &=\sum_{n=0}^{\infty}  \lim_{k \to \infty} G(n,k) \nonumber \\
&= 12 \sum_{n=0}^\infty \frac{1}{4^n} \frac{\left(\frac{1}{2} \right)_n^3}{(1)_n^3}(6n+1)
\lim_{k \to \infty} \frac{1}{k} \frac{\left(\frac{1}{2} \right)_k^2}{(1)_k^2}=\frac{48}{\pi^2},
\end{align}
in which we have used the asymptotic approximation $(k)_n \sim k^n$. The series in (\ref{uniform}) is a Ramanujan series with sum $4/\pi$, see \cite{baruah0}, and we have evaluated the last limit using Stirling's formula. Hence, we have
\[
\sum_{n=0}^{\infty} \frac{\left( \frac{1}{2} \right)_n^3 \left(1+\frac{k}{2} \right)_n \left( \frac{1}{2}+\frac{k}{2} \right)_n }{(1)_n^3 (1+k)_n^2} \frac{\left( \frac{1}{2} \right)_k^2}{(1)_k^2} \, \frac{12k(8n^2+6kn+2n+k)}{2n+k}=\frac{48}{\pi^2}.
\]
For example, taking $k=1$, we obtain a formula that Maple can evaluate, namely
\begin{equation}
\sum_{n=0}^{\infty} \frac{ \left( \frac{1}{2} \right)_n^4}{(1)_n^4} \frac{8n^2+8n+1}{(n+1)^2}=\frac{16}{\pi^2},
\end{equation}
which is an example of series which converge slowly to the constant $1/\pi^2$. To obtain more interesting series,
we replace $k$ with $k+n$ in $F(n,k)$. Then, we have the new function
\begin{equation}\label{segF}
F(n,k)=U(n,k)  \frac{96n^3}{3n+k},
\end{equation}
where
\[
U(n,k)=\left(\frac{27}{64} \right)^n \frac{\left( \frac{1}{2} \right)_n^3 \left( \frac{1}{2}+k \right)_n^2 \left( 1+\frac{k}{3} \right)_n
\left(\frac{1}{3}+\frac{k}{3} \right)_n \left( \frac{2}{3}+\frac{k}{3} \right)_n}
{(1)_n^3 (1+k)_n \left(1+\frac{k}{2} \right)_n^2 \left( \frac{1}{2}+\frac{k}{2} \right)^2}
\frac{\left( \frac{1}{2} \right)_k^2}{(1)_k^2}.
\]
Package EKHAD gives the companion
\begin{equation}\label{segG}
G(n,k)=U(n,k) \frac{n(2n+1)^2(74n^2+27n+3)+kP(n,k)}{(n+\frac{k}{3})(2n+k+1)^2},
\end{equation}
where
\begin{align}
P(n,k) &=(2n+1)(296n^3+164n^2+26n+1) \nonumber \\ &+(480n^3+360n^2+78n+5)k \nonumber \\ &+(176n^2+80n+8)k^2 \nonumber \\ &+(24n+4)k^3.
\nonumber
\end{align}
If we observe the steps in (\ref{uniform}), we see that again we have
\[ \sum_{n=0}^{\infty} G(n,k)=\frac{48}{\pi^2}. \]
Finally, taking $k=0$, we obtain
\begin{equation}
\sum_{n=0}^{\infty} \left( \frac{3}{4} \right)^{3n} \frac{\left(\frac{1}{2} \right)_n^3 \left(\frac{1}{3} \right)_n
\left(\frac{2}{3} \right)_n}{(1)_n^5} (74n^2+27n+3)=\frac{48}{\pi^2}.
\end{equation}
Although the convergence of this series is not very fast, it seems to us very interesting. The reason is that it is a new formula which belongs to a family of series for $1/\pi^2$ discovered by the author. See \cite{guilleraAAMwz}, \cite{guilleraEMpi2}, \cite{guilleraRJgen}, \cite{guilleratesis} and \cite{baruah0}, \cite{borwein}, \cite{zudilin}. Until now the unique existing proofs, and only for some of these series, are based on the WZ-method. However, it would be a major achievement to find a modular-like theory which can  explain all these kind of formulas; see \cite{chen}, \cite{zudilin-modlike} and \cite{guillera-modlike}.

\section{An Hypergeometric identity}
If $F(n,k)$ and $G(n,k)$ is a WZ-pair then obviously $F_x(n,k)=F(n+x,k)$ and $G_x(n,k)=G_(n+x,k)$ is also a WZ-pair
for every value of $x$. If the last limit in (\ref{sumas-wz}) is equal to zero then, if we repeat the proof in \cite{amdeberhan} we see that
\begin{align}
\sum_{n=0}^{\infty} G_x(n ,0) &=\sum_{n=0}^{\infty} G_x(n,1)+ F_x(0,0) = \sum_{n=0}^{\infty} G_x(n,2)+ F_x(0,1)+F_x(0,0) \nonumber \\
&=\sum_{n=0}^{\infty} G_x(n,3)+\sum_{k=0}^{2} F_x(0,k)=\sum_{n=0}^{\infty} G_x(n,4)+\sum_{k=0}^{3} F_x(0,k) =\cdots. \nonumber
\end{align}
Therefore, as in \cite{amdeberhan}, we arrive to
\[
\sum_{n=0}^{\infty} G_x(n ,0)=\lim_{k \to \infty} \sum_{n=0}^{\infty} G_x(n,k)+\sum_{k=0}^{\infty} F_x(0,k).
\]
This is the formula we used to obtain the formulas in \cite{guillera-hyperiden}. Observe that for $x=0$ the last sum is zero. If we now apply the formula to the WZ-pair of functions (\ref{segF}) and (\ref{segG}), we obtain the following hypergeometric identity:
\begin{multline}
\frac{1}{48}\sum_{n=0}^{\infty} \left( \frac{27}{64} \right)^n \frac{\left(\frac{1}{2} \right)_{n+x}^3 \left(\frac{1}{3} \right)_{n+x}
\left(\frac{2}{3} \right)_{n+x}}{(1)_{n+x}^5} (74(n+x)^2+27(n+x)+3)\nonumber \\
=\frac{1}{4 \pi}\sum_{n=0}^{\infty} \left( \frac{1}{4} \right)^{n+x} \frac{\left( \frac{1}{2} \right)_{n+x}^3}{(1)_{n+x}^3} (6(n+x)+1) \\ +2x^3 \left( \frac{27}{64} \right)^x \frac{ \left(\frac{1}{2} \right)_{x}^3 \left(\frac{1}{3} \right)_{x} \left(\frac{2}{3} \right)_x}{(1)_{x}^5}\sum_{k=0}^{\infty} \frac{\left( \frac{1}{2}+x \right)_k^2
(1+3x)_k}{(1+2x)_k^2 (1+x)_k} \, \frac{1}{k+3x}. \nonumber
\end{multline}
or equivalently
\begin{multline}
\frac{1}{48}\sum_{n=0}^{\infty} \left( \frac{27}{64} \right)^n \frac{\left(\frac{1}{2}+x \right)_n^3 \left(\frac{1}{3}+x \right)_n
\left(\frac{2}{3}+x \right)_n}{(1+x)_n^5} (74(n+x)^2+27(n+x)+3) \nonumber \\
=\frac{2x}{\pi} \left( \frac{16}{27} \right)^x \, \frac{(1)_x^2}{\left( \frac{1}{3} \right)_x \left( \frac{2}{3} \right)_x}
\sum_{n=0}^{\infty} \frac{\left( \frac{1}{2} \right)_n^2}{(1+x)_n^2}+
2x^3 \sum_{n=0}^{\infty} \frac{\left( \frac{1}{2}+x \right)_n^2 (1+3x)_n}{(1+2x)_n^2 (1+x)_n} \, \frac{1}{n+3x}. \nonumber
\end{multline}
where we have used \cite[Iden. 1]{guillera-hyperiden}. Taking $x=1/2$, we get
\begin{equation}
\sum_{n=0}^{\infty} \left( \frac{27}{64} \right)^n \frac{(1)_n^3 \left(\frac{1}{6}\right)_n
\left(\frac{5}{6}\right)_n}{\left(\frac{1}{2}\right)_n^5} \frac{(74n^2+101n+35)(6n+1)}{(2n+1)^5}=\frac{16\pi^2}{3},
\end{equation}
which is a new formula for $\pi^2$.

\enddocument